\documentclass{amsart}

\usepackage{amsmath,amssymb}
\usepackage[all]{xy}

\newtheorem{theorem}{Theorem}
\newcommand{\bt}{\begin{theorem}}
\newcommand{\et}{\end{theorem}}
\newtheorem{lemma}{Lemma}
\newcommand{\bl}{\begin{lemma}}
\newcommand{\el}{\end{lemma}}
\newtheorem{corollary}{Corollary}
\newcommand{\bc}{\begin{corollary}}
\newcommand{\ec}{\end{corollary}}
\newcommand{\beq}{\begin{equation}}
\newcommand{\eeq}{\end{equation}}
\newcommand{\benum}{\begin{enumerate}}
\newcommand{\eenum}{\end{enumerate}}
\newcommand{\N}{\ensuremath{ \mathbf N }}

\newcommand{\R}{\ensuremath{\mathbf R}}

\newcommand{\bmat}{\left(\begin{matrix}}
\newcommand{\emat}{\end{matrix}\right)}

\DeclareMathOperator{\qqand}{\qquad\text{and}\qquad}

\title{Pairs of matrices in $GL_2(\R_{\geq 0})$ that freely generate}
\author{Melvyn B. Nathanson}
\address{Department of Mathematics\\
Lehman College (CUNY)\\Bronx, NY 10468}
\email{melvyn.nathanson@lehman.cuny.edu}

\date{\today}

\begin{document}

\begin{abstract}
An elementary proof that certain pairs of $2\times 2$ matrices 
with nonnegative real coordinates generate free monoids.
\end{abstract}

\maketitle

A monoid is a semigroup with an identity.
Let $GL_2(\R_{\geq 0})$ denote the multiplicative monoid of $2 \times 2$ 
matrices with nonzero determinant and with coordinates in 
the set $\R_{\geq 0}$ of nonnegative real numbers.  
To every matrix 
\[
X =  \bmat x_{1,1} & x_{1,2} \\ x_{2,1} & x_{2,2} \emat \in GL_2(\R_{\geq 0})
\]
we associate the  linear fractional transformation
\[
X(t) = \frac{x_{1,1}t + x_{1,2} }{x_{2,1} t+ x_{2,2}} 
\]
defined on the set of positive real numbers $t$.  
This is a monoid isomorphism from $GL_2(\R_{\geq 0})$ 
to the monoid of linear fractional transformations with 
nonnegative real coordinates, nonzero determinant, and the binary operation 
of composition of functions.   Without loss (or gain) of generality, we can use the language of matrices or the language of functions.

It is important to note that the nonzero determinant of $X$ 
and the nonnegativity of the coordinates of $X$ imply that
if $t > 0$, then $X(t) > 0$.

The monoid $M(A,B)$ generated by a pair of matrices 
$\{A,B\}$ in $GL_2(\R_{\geq 0})$
consists of all matrices that can be represented as 
products of nonnegative powers of $A$ and $B$.
The matrices $A$ and $B$ \emph{freely generate} this monoid if every matrix in $M(A,B)$ has a unique representation as a product of powers of $A$ and $B$.

Consider  the matrices 
$L_1 = \bmat 1 & 0 \\ 1 & 1 \emat$ and 
$R_1 = \bmat 1 & 1 \\ 0 & 1 \emat$.
The Calkin-Wilf tree~\cite{calk-wilf00,nath14f} 
is the directed graph whose vertices are 
the positive rational numbers and which is 
constructed inductively from the root vertex 1 by the generation rule 
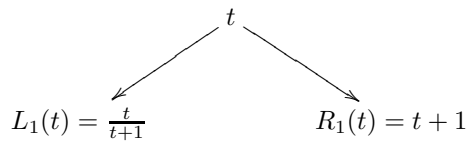
\[
\xymatrix{
 & t \ar[dl] \ar[dr]&  \\
L_1(t) = \frac{t}{t+1}& & R_1(t) = t+1 
}
\]
The fact that the Calkin-Wilf graph is a tree is equivalent to the well-known folk 
theorem that the matrices $L_1$ and $R_1$ freely generate 
the monoid $SL_2(\N_0)$ of $2 \times 2$ 
matrices with  determinant 1 and nonnegative integral coordinates.  

It is a standard application of the ping-pong lemma 
that for every pair $(u,v)$ of integers with $u \geq 2$ and $v \geq 2$, 
the matrices $L_u = \bmat 1 & 0 \\ u & 1 \emat$ and 
$R_v = \bmat 1 & v \\ 0 & 1 \emat$ generate a free group of rank 2 
(Lyndon and Schupp~\cite[pp. 167--168]{lynd-schu77}).  
The case $u=v=2$ is Sanov's theorem~\cite{sano47a}.

In this note we give a simple proof that certain pairs of matrices 
in $GL_2(\R_{\geq 0})$ freely generate a monoid.  
The examples above are special cases of this result.

\bt            \label{CWFree:theorem}
Let $A = \bmat a_{1,1} & a_{1,2} \\ a_{2,1} & a_{2,2} \emat $ 
and  $B = \bmat b_{1,1} & b_{1,2} \\ b_{2,1} & b_{2,2} \emat $ 
be matrices in $GL_2(\R_{\geq 0})$.  
If 
\beq                     \label{CWFree:ineqA}
a_{1,1} \leq a_{2,1} \qqand  a_{1,2} \leq a_{2,2}
\eeq
and if 
\beq                     \label{CWFree:ineqB}
b_{1,1} \geq b_{2,1}  \qqand  b_{1,2} \geq b_{2,2}
\eeq
then the submodule of $GL_2(\R_{\geq 0})$ generated by $A$ and $B$ is free, 
and $\{A,B\}$ is a free set of generators  for this module.  
\et

\begin{proof}
We associate to the matrices $A$ and $B$ the linear fractional transformations
\[
A(t) = \frac{a_{1,1}t + a_{1,2} }{a_{2,1} t+ a_{2,2}}
\qqand
B(t) = \frac{b_{1,1}t + b_{1,2} }{b_{2,1} t+ b_{2,2}}.
\]
Let $t > 0$.  
Because $\det(A) = a_{1,1}a_{2,2} - a_{1,2}a_{2,1} \neq 0$, 
we have $a_{2,1}t + a_{2,2} > 0$.
Inequalities~\eqref{CWFree:ineqA} imply that 
\[
(a_{1,1} - a_{2,1})t \leq 0 \leq a_{2,2} - a_{1,2}.
\]
If $(a_{1,1} - a_{2,1})t = a_{2,2} - a_{1,2}$, then 
$a_{1,1} = a_{2,1}$ and $a_{2,2} = a_{1,2}$, 
and so $\det(A) = 0$, 
which is absurd.  
Therefore, $(a_{1,1} - a_{2,1})t < a_{2,2} - a_{1,2}$ 
or, equivalently, $0 < A(t) < 1$.  
Similarly,  inequalities~\eqref{CWFree:ineqB} imply that if $t' > 0$, 
then $(b_{1,1} - b_{2,1})t' > b_{2,2} - b_{1,2}$ and so 
$B(t') > 1$ .  
Thus, for all $t , t'> 0$, we have  
\beq                     \label{CWFree:ineqABt}
0 < A(t) < 1 < B(t').
\eeq

If $A$ and $B$ do not freely generate a monoid, then there exist distinct sequences 
$(X_1, X_2, \ldots, X_k)$ and $(Y_1, Y_2, \ldots, Y_{\ell})$ with $X_i \in \{A,B\}$ 
for $i=1,\ldots, k$ and $Y_j \in \{A,B\}$  for $j=1,\ldots, \ell$ such that  
\beq                     \label{CWFree:XY}
X_1 X_2 \cdots X_k = Y_1 Y_2 \cdots Y_{\ell}.
\eeq
 Choose the smallest positive integer $k$ for which a relation 
 of the form~\eqref{CWFree:XY} exists.  If $X_1 = Y_1$, then 
 $X_2 \cdots X_k = Y_2 \cdots Y_{\ell}$, which contradicts the minimality of $k$.
 Therefore, $X_1 \neq Y_1$.  
 
Suppose that $X_1 = A$ and $Y_1 = B$. 
Applying the matrices as linear fractional transformations, 
we obtain $t = X_2 \cdots X_k(1) > 0$ and $t' = Y_2 \cdots Y_{\ell}(1) > 0$.
Identity~\eqref{CWFree:XY} implies that  
\[
A(t) = X_1X_2 \cdots X_k(1) = Y_1 Y_2 \cdots Y_{\ell}(1) = B(t').
\]
This is absurd because it contradicts inequality~\eqref{CWFree:ineqABt}.  
The case $X_1 = B$ and $Y_1 = A$ is similar. 
This completes the proof.
\end{proof}

\def\cprime{$'$} \def\cprime{$'$} \def\cprime{$'$} \def\cprime{$'$}
\providecommand{\bysame}{\leavevmode\hbox to3em{\hrulefill}\thinspace}
\providecommand{\MR}{\relax\ifhmode\unskip\space\fi MR }
\providecommand{\MRhref}[2]{%
  \href{http://www.ams.org/mathscinet-getitem?mr=#1}{#2}
}
\providecommand{\href}[2]{#2}

\end{document}